\documentclass[a4paper,reqno]{amsart}%
\usepackage{amssymb,latexsym}
\usepackage{amsmath}
\usepackage{amsfonts}
\usepackage{amssymb}
\usepackage{graphicx}%
\theoremstyle{plain}

\theoremstyle{definition}

\theoremstyle{remark}

\numberwithin{equation}{section}
\numberwithin{theorem}{section}

\begin{document}
\title[Associative Spectra]{Associative Spectra of Binary Operations}
\author{B\'{e}la Cs\'{a}k\'{a}ny}
\address[B\'{e}la Cs\'{a}k\'{a}ny]{Bolyai Institute, University of Szeged, Aradi v\'{e}rtan\'{u}k tere 1, H-6720
Szeged, Hungary}
\email{csakany@math.u-szeged.hu}
\author{Tam\'{a}s Waldhauser}
\address[Tam\'{a}s Waldhauser]{Bolyai Institute, University of Szeged, Aradi v\'{e}rtan\'{u}k tere 1, H-6720
Szeged, Hungary}
\email{twaldha@math.u-szeged.hu}
\thanks{This research was supported by Hungarian National Foundation for Scientific
Research (OTKA) grants no. T022867 and T026243.}
\dedicatory{To Ivo Rosenberg on his 65th birthday}

\begin{abstract}
The distance of a binary operation from being associative can be
\textquotedblleft measured\textquotedblright\ by its \textit{associative
spectrum}, an appropriate sequence of positive integers. Particular instances
and general properties of associative spectra are studied.

\end{abstract}
\maketitle

\section{Introduction}

Let $n$ be a positive integer. We call a string consisting of symbols $x,(,$
and $)$ a \textit{bracketing of size }$n$ if it contains $n$ symbols
\textquotedblleft$x$\textquotedblright, and $n-1$ symbols \textquotedblleft%
$($\textquotedblright\ (left parentheses) as well as \textquotedblleft%
$)$\textquotedblright\ (right parentheses) so that they are properly placed to
determine a product of $n$ factors $x$ (see, \textit{e.g.} [1,15]). More formally,

\smallskip

\noindent1. $x$ is the unique bracketing of size 1,

\noindent2. the bracketings of size $n$ are exactly the strings of form $(PQ)$
where $P$ and $Q$ are bracketings of size $k$ resp. $l$ with $k+l=n$.

\smallskip

\textit{E.g.} $(xx)$ is the only bracketing of size 2, and $((x(xx))(xx))$ is
a bracketing of size 5. Note that we always use an outermost pair of
parentheses whenever $n>1$, in contrary to the everyday usage of parentheses.
We shall denote bracketings by capital letters, and $|B|$ stands for the size
of $B$.

Bracketings are, in fact, the elements of the free groupoid\footnote{{Here by
a groupoid we mean a nonempty set with one binary operation.}} with one free
generator $x$ (\textit{cf.} [1], p. 133), or, equivalently, they are the unary
groupoid terms. The corresponding unary term operations on special groupoids
were investigated by several authors (see, \textit{e.g.} [5,7]). In any
bracketing of size $n$ we can indicate the position of symbols $x$ by
subscripts $1,\ldots,n$, \textit{e.g.} $(x_{1}x_{2}),((x_{1}(x_{2}%
x_{3}))(x_{4}x_{5}))$. Thus, a bracketing of size $n$ provides also an element
of the free groupoid with free generators $x_{1},\ldots,x_{n}$, \textit{i.e.},
an $n$-ary groupoid term (although, of course, not all $n$-ary groupoid terms
originate from bracketings in such a way). Here we always study bracketings
considered as $n$-ary groupoid terms, even if in some cases we omit the
subscripts $1,\ldots,n$. On every groupoid $G$, these terms give rise to
$n$-ary term operations. We call them \textit{regular} $n$\textit{-ary
operations of} $G$ (or, \textit{regular over the operation of} $G$), and, in
concrete cases, \textit{operations induced by given bracketings}. For notation
of the regular operation induced by the bracketing $B,P,Q,$ \textit{etc.} we
use the corresponding lowercase letters $b,p,q,$ \textit{etc}.

If $G$ is associative, then by the generalized associative law there is
exactly one regular $n$-ary operation for each $n$. In the general case, we
have a sequence
\[
{s_{{}_{G}}}(1),{s_{{}_{G}}}(2),\ldots,{s_{{}_{G}}}(n),\ldots
\]
of positive integers with $s_{{}_{G}}(n)$ denoting the number of distinct
$n$-ary regular operations of $G$. \textit{E.g.}, $s_{{}_{G}}(1)=s_{{}_{G}%
}(2)=1 $ for every groupoid $G$, and $s_{{}_{G}}(3)=2$ if and only if $G$ is
nonassociative, as then the two possible bracketings of size 3, $(x_{1}%
(x_{2}x_{3}))$ and $((x_{1}x_{2})x_{3})$ induce different ternary term operations.

The sequence
\[
\{s_{{}_{G}}(n)\}=(s_{{}_{G}}(1),s_{{}_{G}}(2),\ldots,s_{{}_{G}}(n),\ldots)
\]
measures, in some sense, the distance of $G$ from associativity: the smaller
its entries are, the closer the operation of $G$ is to being associative.
Hence we call this sequence \textit{the associative spectrum of} $G$ (or,
\textit{of the operation of} $G$). Instead of $s_{{}_{G}}(n)$ we write $s(n)$
if this cannot cause misunderstanding. Usually we also omit $s(1)$ and $s(2)$,
bearing no information about $G$.

In this paper we study the introduced notion from several points of view. The
next section contains some well-known facts, simple observations, and
auxiliary results on bracketings and associative spectra; there and later, the
routine inductive proofs will often be omitted. Most frequently we use
induction on size; we leave out the words \textquotedblleft on
size\textquotedblright\ in these cases. The third section contains samples of
determining associative spectra of some familiar nonassociative operations.
The problem of characterizing all associative spectra of operations on a set
with a given power seems to be hard. However, the case of the two-element set
is, as it might be expected, easy (Section 4), and a lot of three-element
groupoids are accessible (Section 5). In the final section we present some
facts on the general behavior of associative spectra, and formulate several problems.

Further on, we write simply \textit{spectrum} for associative spectrum.

\section{Properties of bracketings and spectra}

For any bracketing $B$ of size $n(>1)$, we can \textit{pair} its left and
right parentheses in a natural way ([9,15]). Induction shows that we can
always choose a consecutive quadruple $(xx)$ in $B$; its left and right
parentheses will be associated to form a pair. Replacing then $(xx)$ with $x$
we obtain a bracketing $B^{\prime}$ of size $n-1$, for which the preceding
process can be repeated until no unpaired parentheses remain. This way of
forming pairs involves that any pair together with the symbols between them is
also a bracketing. It is called a \textit{subbracketing of} $B$;
\textit{e.g.}, if $B=(PQ)$, then $P$ and $Q$ are subbracketings of $B$, as
outermost parentheses of any bracketing are paired. We call $P$ and $Q$ the
(\textit{left} resp. \textit{right}) \textit{factors} of $B$. The symbols $x$
are considered as subbracketings of size 1, too. Observe that pairing is
unique, and if a parenthesis lies between a pair then its associate also lies
between them. Hence the representation of bracketings of size $>1$ in form
$(PQ)$ is unique, too.

Substituting $x$ for one or several disjoint subbracketings in $B$ we obtain
\textit{a quotient bracketing} of $B$. \textit{E.g.} $(x(xx))$ and
$((xx)(xx))$ are (disjoint) subbracketings of $B=(((x(xx))x)((xx)(xx)))$, and
replacing them with $x$ provides the quotient bracketing $((xx)x)$ of $B$. A
bracketing is a \textit{nest} if it is either of size 1 (\textit{a trivial
nest}) or one out of its factors is $x$, and the other one is a nest ([5,7]).
\textit{E.g.}, all bracketings of size 4 save $(xx)(xx)$ are nests. Given a
bracketing $B$, there are subbracketings of $B$ which are nests; in
particular, each $x_{i}$ is contained in a unique maximal nest. We call these
maximal nests simply \textit{the nests of} $B$. A nontrivial nest has a unique
subbracketing of form $(x_{i}x_{i+1})$; we say that $x_{i},x_{i+1}$ are
\textit{the eggs of the nest}.

The \textit{Catalan numbers} $C_{n}$ are defined recursively by

\smallskip

\noindent$\left(  1\right)  $ $C_{0}=1$,

\noindent$\left(  2\right)  $ $C_{n}=C_{0}C_{n-1}+C_{1}C_{n-2}+\cdots
+C_{n-2}C_{1}+C_{n-1}C_{0}\quad(n>0)$,

\smallskip

\noindent or, equivalently, by the formula%
\[
C_{n}={\frac{1}{{n+1}}}{\binom{{2n}}{n}}.
\]
Compare (1) and (2) with the formal definition of bracketings in the
introduction, and take into account the unicity of the representation of
bracketings in form $(PQ)$. Then the following standard result follows:

\bigskip

\noindent2.1.\textit{\ The number of bracketings of size }$n$\textit{\ equals
}$C_{n-1}$ (see, \textit{e.g.} [8]).

\smallskip

\noindent Hence we infer:

\bigskip

\noindent2.2.\textit{\ For any spectrum }$\{s(n)\}$\textit{,}%
\[
1\leq s(n)\leq C_{n-1}%
\]
\textit{holds for every }$n$\textit{.}

\smallskip

If $s_{{}_{G}}(n)=C_{n-1}$ for every $n$, then the groupoid $G$ and its
operation are said to be \textit{Catalan}. \textit{E.g.}, free groupoids are
Catalan. Another inequality also follows from the definition of bracketings:

\bigskip

\noindent2.3.\textit{\ For any spectrum }$\{s(n)\}$\textit{, }%
\[
s(n)\leq s(1)s(n-1)+s(2)s(n-2)+\cdots+s(n-2)s(2)+s(n-1)s(1)
\]
\textit{holds for every }$n(\geq2)$\textit{.}

\smallskip

Hence if $s_{{}_{G}}(n_{0})<C_{{n_{0}}-1}$ then $s_{{}_{G}}(n)<C_{n-1}$ for
every $n>n_{0}$. The following trivial observations are useful, too:

\bigskip

\noindent2.4.\textit{\ If the groupoids }$G$\textit{\ and }$H$\textit{\ are
isomorphic or antiisomorphic, then their spectra coincide.}

\bigskip

\noindent2.5.\textit{\ If the groupoid }$H$\textit{\ is a subgroupoid or a
factorgroupoid of }$G$\textit{, then }%
\[
s_{{}_{H}}(n)\leq s_{{}_{G}}(n)
\]
\textit{holds for every }$n$\textit{.}

\smallskip

By 2.5, in order to show that $G$ is Catalan it is sufficient to find a
Catalan subgroupoid or factorgroupoid of $G$. The next fact goes back to
\L ukasiewicz (for a proof, see [3], Ch. 3.2, or [11], Exercise 1.38):

\bigskip

\noindent2.6.\textit{\ Bracketings are uniquely determined by the places of
their right (or left) parentheses between the symbols }$x_{1},\ldots,x_{n}%
$\textit{.}

\smallskip

Next we introduce sequences of nonnegative integers which arise naturally from
bracketings, and also contain full information on them. Consider the free
monoid $F_{2}$ with unit element $e$, generated by symbols 0 and 1. A subset
$M$ of $F_{2}$ is \textit{prefix-free} if no word in $M$ is a prefix
(\textit{i.e.}, a left segment) of another word in $M$. There exist finite
maximal prefix-free sets (\textit{FMPF-sets} in short) in $F_{2}$,
\textit{e.g.}, the set containing the empty word $e$ only, the sets
$\{0,1\},\{00,010,011,10,11\},$ \textit{etc}. Assign to each bracketing an
ordered sequence of words in $F_{2}$ inductively by the rule:

\smallskip

\noindent(a) $x\mapsto(e)$,

\noindent(b) if$\hfill P\,\mapsto\,(w_{1},\ldots,w_{k})\hfill$and$\hfill
Q\,\mapsto\,(w_{k+1},\ldots,w_{k+l})\hfill$then$\hfill(PQ)\mapsto
(0w_{1},\ldots,0w_{k},$

$~1w_{k+1},\ldots,1w_{k+l}).$

\smallskip

It is a routine to check that, in this way, a unique, lexicographically listed
FMPF-set of $n$ words is assigned to every bracketing of size $n$. Now we can
use the defining properties (1),(2) of Catalan numbers to show that the number
of distinct FMPF-sets of $n$ elements equals $C_{n-1}$. Therefore, (a) and (b)
provide a 1-1 correspondence between bracketings and lexicographically ordered FMPF-sets.

\pagebreak

Consider a bracketing $B$ of size $n$ viewed with subscripts, \textit{i.e.},
as an $n$-ary groupoid term. Let $({w_{1}}(B),\ldots,{w_{n}}(B))$ be the
lexicographically ordered FMPF-set corresponding to $B$. Call the length of
${w_{i}}(B)$ \textit{the depth of} $x_{i}$ \textit{in} $B$, and the number of
$0$'s (resp. of $1$'s) in ${w_{i}}(B)$ \textit{the left depth} (resp.
\textit{the right depth}) \textit{of} $x_{i}$ \textit{in} $B$.

Inspecting (a) and (b) we get the intuitive meaning of depth of $x_{i}$: the
number of pairs of parentheses (or, equivalently, of the subbracketings of
size at least 2\thinspace) containing $x_{i}$. Similarly, \textit{e.g.} the
right depth of $x_{i}$ in $B$ is the number of those subbracketings in which
$x_{i}$ is contained in the right factor. The sequence consisting of the
depths of $x_{1},\ldots,x_{n}$ in $B$ will be called \textit{the depth
sequence of }$B$. \textit{Left} and \textit{right depth sequences of} $B$ are
defined analogously. \textit{E.g.}, the depth sequence of $((x(xx))(xx))$ is
$(2,3,3,2,2)$, and its right depth sequence is $(0,1,2,1,2)$.%

\begin{figure}
[ptb]
\begin{center}
\includegraphics[
width=3.2in
]%
{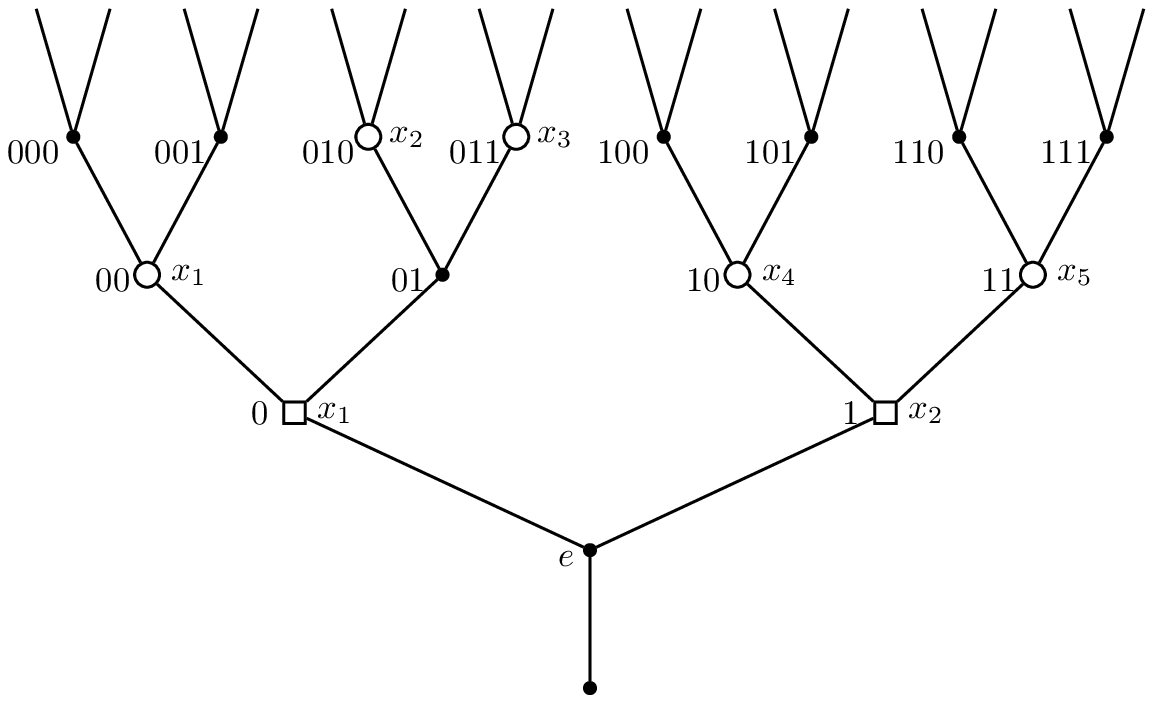}%
\end{center}
\end{figure}

FMPF-sets --- and thus also bracketings --- can be imagined as such minimal
sets of vertices in the infinite binary rooted tree that separate the top of
the tree from its bottom. See the figure where the sets of vertices
corresponding to $(x_{1}x_{2})$ and $(x_{1}(x_{2}x_{3}))(x_{4}x_{5})$ are
marked by squares, resp. circles; correspondence between vertices and binary
strings is indicated, too. In this representation, the depth of $x_{i}$ is the
number of edges in the path $p$ connecting $e$ with $x_{i}$. Similarly, the
left (right) depth is the number of left(right) edges in $p$.

\bigskip

\noindent2.7.\textit{\ Bracketings are uniquely determined by their depth
sequences.}

\smallskip

This is clearly true for bracketings of size $\leq3$. Suppose the bracketings
$(P_{1}Q_{1})$ and $(P_{2}Q_{2})$ of size $n(>3)$ have the same depth sequence
$(d_{1},\ldots,d_{n})$. From the definition, the equality%
\begin{equation}
\sum_{i=1}^{n}{\frac{1}{2^{e_{i}}}}=1 \tag{1}%
\end{equation}
follows for every depth sequence $(e_{1},\ldots,e_{n})$. If $|P_{1}%
|=j,|P_{2}|=k$, then, in view of (a) and (b), the depth sequences of $P_{1}$
and $P_{2}$ are of form $(d_{1}-1,\ldots,d_{j}-1)$ and $(d_{1}-1,\ldots
,d_{k}-1)$, respectively. Therefore,%
\[
\sum_{i=1}^{j}{\frac{1}{2^{d_{i}}}}=\sum_{i=1}^{k}{\frac{1}{2^{d_{i}}}}=1/2.
\]
Hence the sizes of $P_{1}$ and $P_{2}$ are equal. Now the proposition follows
by induction.

\bigskip

\noindent2.8.\textit{\ Bracketings are uniquely determined by their right (or
left) depth sequences.}

\smallskip

Let $B=(PQ)$ be a bracketing with right depth sequence (in short,
\textit{RD-sequence})%
\begin{equation}
(d_{1},\ldots,d_{n}). \tag{2}%
\end{equation}
Then there is a $k$ between $1$ and $n$ such that the RD-sequence of $P$ is
$(d_{1},\ldots,d_{k})$, and that of $Q$ is $(d_{k+1}-1,\ldots,d_{n}-1)$.
Induction shows that always%
\begin{equation}
d_{1}=0,\quad d_{2}=1, \tag{3}%
\end{equation}
and, for $i=1,\ldots,n-1$,%

\begin{equation}
1\leq d_{i+1}\leq d_{i}+1. \tag{4}%
\end{equation}
Call a sequence (2) of nonnegative integers a \textit{zag sequence}
(\textit{cf.} [6], Ch. 1.2, where \textit{zig} is defined) if it has the
properties (3) and (4). We use induction to prove that for any zag sequence
(2) there exists at most one bracketing with RD-sequence (2). This is clearly
true for $n\leq2$. As $(d_{k+1}-1,\ldots,d_{n}-1)$ is a zag sequence, we have
$d_{k+1}=1$, and $d_{j}\geq2$ for $j=k+2,\ldots,n$. It follows that if the
size of the first factor of $B$ is $k$, then the last $1$ in the RD-sequence
of $B$ appears on the $(k+1)$st place. Hence if the RD-sequences of $B=(PQ)$
and $B^{\prime}=(P^{\prime}Q^{\prime})$ are the same, then $|P|=|P^{\prime}|$.
Thus the RD-sequences of $P$ and $P^{\prime}$ coincide, and, by induction,
$P=P^{\prime}$. Similarly we obtain $Q=Q^{\prime}$, completing the proof.

An analogous straightforward induction shows that every zag sequence is the
RD-sequence of some bracketing. Consequently, the number of zag sequences of
length $n$ equals that of the bracketings of size $n$, \textit{i.e.,}
$C_{n-1}$ (\textit{cf.} [14], Ch. 5, Exercise 19(u)).

\section{Examples}

In this section we determine spectra of several common operations. Given a
particular operation, we denote the members of its spectrum by $s(n)$ (without
subscript), and we write $s(n)=f(n)$ to indicate that this equality holds for
$n\geq3$.

\bigskip

\noindent3.1.\textit{\ For the subtraction of numbers, }$s(n)=2^{n-2}%
$\textit{.}

\smallskip

Induction shows that any regular operation $b(x_{1},x_{2},\ldots,x_{n})$ over
the subtraction is of form $x_{1}-x_{2}\pm x_{3}\pm\cdots\pm x_{n}$. It is
enough to prove that actually every possible sequence of the $+$ and $- $
signs occurs. This is true for $n\leq3$; suppose $n>3$, and apply induction.
If $b(x_{1},x_{2},\ldots,x_{n})=x_{1}-x_{2}-\cdots-x_{n}$, then $b$ is induced
by $((\ldots((x_{1}x_{2})x_{3})\ldots)x_{n})$. Otherwise there exists a first
$+$ sign in $f$, say $b(x_{1},x_{2},\ldots,x_{n})=x_{1}-x_{2}-\cdots
-x_{k+1}+x_{k+2}\pm\cdots\pm x_{k+l}~(k+l=n)$. Then $b(x_{1},\ldots
,x_{n})=(x_{1}-x_{2}-\cdots-x_{k})-(x_{k+1}-x_{k+2}\mp\cdots\mp x_{k+l})$, and
this is induced by $B=(PQ)$, where $P=((\ldots((x_{1}x_{2})x_{3})\ldots
x_{k})$, and $Q$ is the bracketing that induces the subtrahend (such a $Q$
exists by induction). In fact, this reasoning is valid for subtraction in
arbitrary Abelian groups except those of exponent 2.

\bigskip

\noindent3.2.\textit{\ The arithmetic mean as a binary operation on numbers is
Catalan.}

\smallskip

We prove that distinct bracketings induce distinct regular operations over the
arithmetic mean. Induction shows that a bracketing $B$ of size $n$
induces\linebreak$b(x_{1},\ldots,x_{n})=\sum_{i=1}^{n}{2^{-d_{i}}x_{i}}$ over
the arithmetic mean, where $d_{i}$ is the depth of $x_{i}$ in $B$. Let
$B^{\prime}\,(\neq B)$ be another bracketing of size $n$ which induces
$b^{\prime}(x_{1},\ldots,x_{n})=\sum_{i=1}^{n}{2^{{-d_{i}}^{\prime}}x_{i}}$.
In virtue of 2.7, there exists a $j~(1\leq j\leq n)$ such that $d_{j}%
\neq{d_{j}}^{\prime}$. Then $b(\delta_{1}^{j},\ldots,\delta_{n}^{j}%
)=2^{-d_{j}}\neq{2^{{-d_{j}}^{\prime}}}=b^{\prime}(\delta_{1}^{j}%
,\ldots,\delta_{n}^{j})$, \textit{i.e.}, $b$ and $b^{\prime}$ are distinct
operations, as required. This holds for an arbitrary set of numbers closed
under arithmetic mean, containing more than one element.

\bigskip

\noindent3.3.\textit{\ The geometric mean and the harmonic mean as binary
operations on positive real numbers are Catalan.}

\smallskip

This follows from 3.2 and 2.4, as the groupoids $(\mathbf{R},\left(
{x+y}\right)  /2)$ and $(\mathbf{R}_{+},\sqrt{xy})$ are isomorphic, as well as
$(\mathbf{R}_{+},\left(  x+y\right)  /2)$ and $({\mathbf{R}_{+}},{2xy/}\left(
x+y\right)  )$.

\bigskip

\noindent3.4.\textit{\ The exponentiation as a binary operation }$(a,b)\mapsto
a^{b}$\textit{\ on numbers is Catalan.}

\smallskip

Let $p_{1},\ldots,p_{n}$ be distinct prime numbers. Consider bracketings
$B,\,B^{\prime}(\neq B)$ and the regular operations $b,b^{\prime}$ they induce
over the exponentiation. We show that $b\neq b^{\prime}$. Making use of the
law $(r^{s})^{t}=r^{st}$, and the usual convention of writing $r^{s^{t}}$
instead of $r^{(s^{t})}$, we can write expressions of form $b(p_{1}%
,\ldots,p_{n})$ without parentheses, \textit{e.g.}, if $B=\left(  \left(
x_{1}\left(  x_{2}x_{3}\right)  \right)  \left(  x_{4}x_{5}\right)  \right)  $
and $p_{i}$ are the first primes, we have $b(2,3,5,7,11)=2^{{3^{5}}{7^{11}}}$.
Here the exponents are at different levels: say, 2 is at the zeroth, 3 and 7
are at the first level, \textit{etc}. The key observation is that the height
of the level of $p_{i}$ in $b$ always equals the right depth of $x_{i}$ in
$B$; this can be verified using induction. As $B\neq B^{\prime}$, by 2.8.
there exists a $j$ such that the right depth of $x_{j}$ in $B$ is different
from that of $x_{j}$ in $B^{\prime}$. Then the fundamental theorem of
arithmetic implies $b(p_{1},\dots,p_{n})\neq b^{\prime}(p_{1},\ldots,p_{n})$.

\bigskip

\noindent3.5.\textit{\ The cross product of vectors is Catalan.}

\smallskip

Consider three pairwise perpendicular unit vectors, their additive inverses,
and the zero vector. They form a groupoid under cross product, and, if we
identify the unit vectors with their negatives, we obtain a four-element
factorgroupoid $C$ with Cayley operation table%

\[
\renewcommand{\arraystretch}{1.4} \setlength{\tabcolsep}{0.2cm}%
\begin{tabular}
[c]{c|c|c|c|c|}%
$\times$ & $0$ & $u$ & $v$ & $w$\\\hline
$0$ & $0$ & $0$ & $0$ & $0$\\\hline
$u$ & $0$ & $0$ & $w$ & $v$\\\hline
$v$ & $0$ & $w$ & $0$ & $u$\\\hline
$w$ & $0$ & $v$ & $u$ & $0$\\\hline
\end{tabular}
\]

\bigskip

In virtue of 2.5., it is enough to prove that this operation is Catalan. Let
$B,B^{\prime},b,b^{\prime}$ be as in 3.4. We shall find nonzero elements
$c_{1},\ldots,c_{n}\in C$ such that $b(c_{1},\ldots,c_{n})=0\neq b^{\prime
}(c_{1},\ldots,c_{n})$. The case $n=3$ is obvious. The general case needs some preparations:

\smallskip

\noindent3.5.1. Let $F$ be a nontrivial nest of size $k$ which induces the
regular operation $f$ on $C$. Given $i~(1\leq i\leq k)$, and $c,d\in C$ with
$d\notin\{0,c\}$, we can choose elements $c_{1},\ldots,c_{i-1},c_{i+1}%
,\ldots,c_{k}\in C$ so that $f(c_{1},\ldots,c_{i-1},c,c_{i+1},\ldots,c_{k})=d$.

This is valid also for any bracketing $B$ and its induced regular operation
$b$ instead of $F$ and $f$. Indeed, apply 3.5.1 to the nest of $B$ containing
$x_{i}$, if this nest is nontrivial, and replace this nest by $x$; while if
$x_{i}$ is a trivial nest, replace the eggs of another nest by $x$. Then, in
both cases, use induction for the quotient bracketing. We remark that this
generalized form of 3.5.1 implies that any regular operation over the cross
product is surjective (\textit{i.e.}, it maps $C^{n}$ onto $C$; in fact, this
is the case for all surjective binary operations, \textit{cf.} 4.2.1).

\smallskip

\noindent3.5.2. If $x_{j},x_{j+1}$ are no eggs of any nest of a bracketing
$B$, we can choose $d_{1},\ldots,d_{j-1},\allowbreak d,d_{j+2},\ldots,d_{k}$
in $C$ such that $f(d_{1},\ldots,d_{j-1},d,d,d_{j+2},\ldots,d_{k})\neq0$.

\smallskip

If $B=(PQ)$ with $|P|=k$ and $j+1\leq k$, then for suitable elements
$d,d_{i}\in C$ by induction we have $p(d_{1},\ldots,d_{j-1},d,d,d_{j+2}%
,\ldots,d_{k})=e\neq0$. Now by 3.5.1 there are $d_{k+1},\ldots,d_{n}\in C$
such that $q(d_{k+1},\ldots,d_{n})=f\neq0,e$. Then $b(d_{1},\ldots
,d,d,\ldots,d_{n})=e\times f\neq0$. The case $k<j$ can be treated in a similar
way. Finally, suppose $k=j$. Let us fix $d\neq0$, and apply 3.5.1 to $P$ and
$Q$ with $i=k$ and $i=k+1$, respectively. Then we have elements $d_{1}%
,\ldots,d_{k-1},d_{k+2},\ldots,d_{n}\in C$ such that $p(d_{1},\ldots
,d_{k-1},d)=e$ and $q(d,d_{k+2},\ldots,d_{n})=f$, where $C=\{0,d,e,f\}$. Thus
$b(d_{1},\ldots,d,d,\ldots,d_{n})=e\times f=d\neq0$, completing the proof of 3.5.2.

In order to prove 3.5., first suppose that there is an $i~(1\leq i\leq n)$
such that $x_{i}$ and $x_{i+1}$ are the eggs of a nest of $B$ as well as of
$B^{\prime}$. Replacing $(x_{i}x_{i+1})$ by $x$ in $B$ and $B^{\prime}$, we
obtain quotient bracketings $B_{1}$ resp. ${B_{1}}^{\prime}$ of size $n-1$
with induced regular operations $b_{1}$ and ${b_{1}}^{\prime}$. By induction,
there exist nonzero elements $e_{1},\ldots,e_{n-1}\in C$ such that ${b_{1}%
}(e_{1},\ldots,e_{i},\ldots,e_{n-1})=0\neq{b_{1}}^{\prime}(e_{1},\ldots
,e_{i},\ldots,e_{n-1})$. Let $e^{\prime},e^{\prime\prime}\in C$ be distinct,
and different from $0$ and $e_{i}$. Then $e^{\prime}\times e^{\prime\prime
}=e_{i}$ , and $b(e_{1},\ldots,e_{i-1},e^{\prime},e^{\prime\prime}%
,e_{i+1},\ldots,e_{n-1})=b_{1}(e_{1},\ldots,e_{i},\ldots,e_{n-1})=0\neq
b^{\prime}(e_{1},\ldots,e_{i-1},\allowbreak e^{\prime},e^{\prime\prime
},e_{i+1},\ldots,e_{n-1})$.

Now suppose that no nests of $B$ and $B^{\prime}$ have a common pair of eggs.
Let $x_{j}$ and $x_{j+1}$ be the eggs of a nest of $B$. Then $b(d_{1}%
,\ldots,d_{j-1},d,d,d_{j+1},\ldots,d_{n})=0$ for any choice of $d_{1}%
,\ldots,d_{j-1},d,d_{j+2},\ldots,d_{n}\in C$. However, as $x_{j}$ and
$x_{j+1}$ are eggs of no nest in $B^{\prime}$, from 3.5.2 it follows that
there is a choice of $d_{1},\ldots,d_{j-1},d_{j+2},\ldots,d_{n}$ such that
$b^{\prime}(d_{1},\ldots,d_{j-1},d,d,d_{j+2},\ldots,d_{n})\neq0$.

\section{Groupoids on two-element sets}

In what follows we consider operations on finite sets. For uniform treatment,
we study groupoids of form $(\mathbf{n},\circ)$, where $\mathbf{n}$ stands for
the set $\{0,1,\ldots,n-1\}$. Each two-element groupoid is isomorphic or
antiisomorphic with $(\mathbf{2},\circ)$, where $x\circ y$ is one of the
following seven Boolean functions:

(1) the constant 1 operation; (2) $x$ (the first projection); (3) $x\wedge y$
(\textit{i.e.}, $\min(x,y)$); (4) $x+y$\thinspace$\operatorname{mod}$ 2; (5)
$x+1$\thinspace$\operatorname{mod}$ 2; (6) $x|y$ (the Sheffer function:
\textquotedblleft neither $x$, nor $y$\textquotedblright); (7) $x\rightarrow
y$ (implication).

Here (1) --- (4) are associative. We determine the spectra of (5) --- (7).

\bigskip

\noindent4.1.\textit{\ For the operation }$x+1\,\operatorname{mod}2$%
,$~s(n)=2$\textit{.}

\smallskip

Indeed, induction shows that for an arbitrary bracketing $B$ of size $n$ and
$\allowbreak c_{1},\ldots,c_{n}\in\mathbf{2}$, $b(c_{1},\ldots,c_{n}%
)=c_{1}+d\,\operatorname{mod}2$, where $d$ is the left depth of $x_{1}$ in $B$.

\bigskip

\noindent4.2.\textit{\ The Sheffer function is Catalan.}

\smallskip

Recall, that $0|0=1$ and $x|y=0$ otherwise. We shall need some preliminaries.

\smallskip

\noindent4.2.1. Regular operations over a surjective operation are surjective
(\textit{i.e.}, they take on all elements of their base sets; the inductive
proof is trivial).

\smallskip

\noindent4.2.2. If the Cayley table of a surjective operation $\circ$ has
neither two identical columns nor two identical rows, then each variable of
any regular operation over $\circ$ is essential.

\smallskip

This is obvious for at most binary regular operations. Let $B=(PQ)$,
$|B|=n\geq3,|P|=k$. Take a variable $x_{i}$ of $b$. We have to prove that
there are elements $c_{1},\ldots,c_{i-1},u,v,c_{i+1},\ldots,c_{n}$ in the base
set $M$ of $\circ$ such that $b(c_{1},\ldots,c_{i-1},u,c_{i+1},\allowbreak
\ldots,c_{n})\neq b(c_{1},\ldots,c_{i-1},v,c_{i+1},\ldots,c_{n})$. Without
loss of generality, suppose $i\leq k$. Then, by induction there exist
$c_{1},\ldots,c_{i-1},u,v,c_{i+1},\ldots,c_{k}\in M$ such that $g=p(c_{1}%
,\ldots,c_{i-1},u,c_{i+1},\ldots,c_{n})\neq p(c_{1},\ldots,c_{i-1}%
,v,c_{i+1},\ldots,c_{n})=h$. The rows of $g$ and $h$ in the Cayley table of
$\circ$ are not identical, \textit{i.e.}, there is a $d\in M$ such that
$g\circ d\neq h\circ d$. Further, by 4.2.1, there are $c_{k+1},\ldots,c_{n}\in
M$ with $q(c_{k+1},\ldots,c_{n})=d$. Then $b(c_{1},\ldots,c_{i-1}%
,u,c_{i+1},\ldots,c_{n})=g\circ d\neq h\circ d=b(c_{1},\ldots,c_{i-1}%
,v,c_{i+1},\ldots,c_{n})$, which was needed.

\smallskip

\noindent4.2.3. If $\circ$ fulfils the conditions of 4.2.2, then regular
operations of distinct arities over $\circ$ cannot be identically equal.

\smallskip

Indeed, otherwise the last variable of the regular operation of greater arity
could not be essential.

We see that 4.2.1---3 apply to the Sheffer function. Let $B_{1},B_{2}$ be
bracketings of size $n\,(\geq3),B_{1}=(P_{1}Q_{1}),B_{2}=(P_{2}Q_{2})$, and
suppose that their induced operations $b_{1}$ and $b_{2}$ coincide. We have to
prove $B_{1}=B_{2}$. This is true for $n=3$, as $(0|0)|1=0\neq1=0|(0|1)$. Let
$n>3$, and assume $k=|P_{1}|\leq|P_{2}|=l$. First we show that, for arbitrary
$c_{1},\ldots,c_{k},\ldots,c_{l}\in\mathbf{2},p_{1}(c_{1},\ldots,c_{k})=0$ if
and only if $p_{2}(c_{1},\ldots,c_{l})=0$. Let $p_{1}(c_{1},\ldots,c_{k})=0$.
By 4.2.1, there exist $c_{k+1},\ldots,c_{n}\in\mathbf{2}$ with $q_{1}%
(c_{k+1},\ldots,c_{n})=0$. Hence it follows
\begin{gather*}
b_{1}(c_{1},\ldots,c_{k},c_{k+1},\dots,c_{n})=p_{1}(c_{1},\ldots
,c_{k})\,|\,q_{1}(c_{k+1},\ldots,c_{n})=1=\\
=b_{2}(c_{1},\ldots,c_{l},c_{l+1},\ldots,c_{n})=p_{2}(c_{1},\ldots
,c_{l})\,|\,q_{2}(c_{l+1},\ldots,c_{n}),
\end{gather*}
implying $p_{2}(c_{1},\ldots,c_{l})=0$. This reasoning is valid in the
opposite direction, too, showing that $p_{1}$ identically equals $p_{2}$. Now
from 4.2.3 we infer $k=l$ and, by induction, $P_{1}=P_{2}$. It remains to
establish $Q_{1}=Q_{2}$. Let, once more, $p_{1}(c_{1},\ldots,c_{k})=0$. If
$Q_{1}\neq Q_{2}$, then, again by induction, there are $c_{k+1},\ldots
,c_{n}\in\mathbf{2}$ such that $q_{1}(c_{k+1},\ldots,c_{n})\neq q_{2}%
(c_{k+1},\ldots,c_{n})$. Then
\[
b_{1}(c_{1},\ldots,c_{n})=0\,|\,q_{1}(c_{k+1},\ldots,c_{n})\neq0\,|\,q_{2}%
(c_{k+1},\ldots,c_{n})=b_{2}(c_{1},\ldots,c_{n}),
\]
a contradiction. Thus $Q_{1}=Q_{2}$, as required.

\bigskip

\noindent4.3.\textit{\ Implication is Catalan.}

\smallskip

In view of 2.4., instead of implication we can consider the operation $x\ast
y$, defined by $0\ast1=1$ and $x\ast y=0$ otherwise, as $(\mathbf{2}%
,\rightarrow)$ and $(\mathbf{2},\ast)$ are isomorphic. For $\ast$, the proof
of 4.2. can be literally adapted.

\section{Groupoids on three-element sets}

There are 3330 essentially distinct three-element groupoids in the sense that
each three-element groupoid is isomorphic with exactly one of them (see the
Siena Catalog [2], in which code numbers from $1$ to $3330$ are given to each
of these representatives), therefore a plain survey of their spectra such as
in the two-element case seems to be impossible. In this section we determine
the spectra of all groupoids on $\mathbf{3}$ with minimal clones of term
operations, and give examples for further spectra.

There exist 12 essentially distinct groupoids on $\mathbf{3}$ with minimal
clones, and each of them is idempotent (see [4]). The operations of an
\textit{idempotent} groupoid on $\mathbf{3}$ may be encoded by the numbers
$0,1,\ldots,728$ in the following transparent way: let the code of $\circ$ be%
\[
(0\circ1)\cdot3^{5}+(0\circ2)\cdot3^{4}+(1\circ0)\cdot3^{3}+(1\circ
2)\cdot3^{2}+(2\circ0)\cdot3+(2\circ1)
\]
(see the examples below). The operations of the groupoids on $\mathbf{3}$ with
minimal clones are (or, more exactly, may be chosen as)
$0,8,10,11,16,17,26,33,35,68,178,624$ (their codes in the Siena Catalog are
$80,102,105,106,122,125,147,267,271,356,\allowbreak1108,2346$ respectively).
It is easy to check that $0,8,10,11$ and $26$ are associative. Here we display
the Cayley tables of the remaining seven operations:%

\begin{gather*}
\renewcommand{\arraystretch}{1.4}\setlength{\tabcolsep}{0.2cm}%
\begin{tabular}
[c]{cccc}%
$%
\begin{tabular}
[c]{|ccc}\hline
$0$ & $0$ & $0$\\
$0$ & $1$ & $1$\\
$2$ & $1$ & $2$%
\end{tabular}
$ & $%
\begin{tabular}
[c]{|ccc}\hline
$0$ & $0$ & $0$\\
$0$ & $1$ & $1$\\
$2$ & $2$ & $2$%
\end{tabular}
$ & $%
\begin{tabular}
[c]{|ccc}\hline
$0$ & $0$ & $0$\\
$1$ & $1$ & $0$\\
$2$ & $0$ & $2$%
\end{tabular}
$ & $%
\begin{tabular}
[c]{|ccc}\hline
$0$ & $0$ & $0$\\
$1$ & $1$ & $0$\\
$2$ & $2$ & $2$%
\end{tabular}
$\\
$16$ & $17$ & $33$ & $35$%
\end{tabular}
\\
\renewcommand{\arraystretch}{1.4}\setlength{\tabcolsep}{0.2cm}%
\begin{tabular}
[c]{ccc}%
$%
\begin{tabular}
[c]{|ccc}\hline
$0$ & $0$ & $0$\\
$2$ & $1$ & $1$\\
$1$ & $2$ & $2$%
\end{tabular}
$ & $%
\begin{tabular}
[c]{|ccc}\hline
$0$ & $0$ & $2$\\
$0$ & $1$ & $1$\\
$2$ & $1$ & $2$%
\end{tabular}
$ & $%
\begin{tabular}
[c]{|ccc}\hline
$0$ & $2$ & $1$\\
$2$ & $1$ & $0$\\
$1$ & $0$ & $2$%
\end{tabular}
$\\
$68$ & $178$ & $624$%
\end{tabular}
\end{gather*}
As we apply three different approaches, we parcel our task into three parts.

\bigskip

\noindent5.1.\textit{\ The operations }16, 17\textit{\ and }178\textit{\ are
Catalan.}

\smallskip

Observe that $\mathbf{3}$ with each of the operations 16, 17 and 178 is a
groupoid in which $\{0,1\}$ is a subgroupoid with two-sided zero element 0,
while $\{1,2\}$ and $\{2,0\}$ are subgroupoids with left zero elements 1 and
2, respectively. Here and in what follows, the just considered operations will
be denoted by circle.

Let $B_{i}=(P_{i}Q_{i})~(i=1,2)$ be distinct bracketings of size $n\,(\geq3)
$. For $n=3$, $1\circ(2\circ0)=1\circ2=1\neq0=1\circ0=(1\circ2)\circ0$,
\textit{i.e.}, $b_{1}\neq b_{2}$. To prove the same for $n>3$, first suppose
$|P_{1}|=k<l=|P_{2}|$. Then%
\begin{align*}
b_{1}(1,\ldots,1,2,\ldots,2,0,\ldots,0)  &  =p_{1}(1,\ldots,1)\circ
q_{1}(2,\ldots,2,0,\ldots,0)=1\circ2=1,\\
b_{2}(1,\ldots,1,2,\ldots,2,0,\ldots,0)  &  =p_{2}(1,\ldots,1,2,\ldots,2)\circ
q_{2}(0,\ldots,0)=1\circ0=0.
\end{align*}

Thus, we can assume $|P_{1}|=|P_{2}|=k$. If $P_{1}\neq P_{2}$, by induction
there exist elements $c_{1},\ldots,c_{k}\in\mathbf{3}$ with $g_{1}=p_{1}%
(c_{1},\ldots,c_{k})\neq p_{2}(c_{1},\ldots,c_{k})=g_{2}$. Let $d$ be the
element of $\mathbf{3}$ that is different from $g_{1}$ and $g_{2}$. Then
$g_{1}\circ d\neq g_{2}\circ d$ (see the Cayley tables), and hence
$b_{1}(c_{1},\ldots,c_{k},d,\ldots,d)=g_{1}\circ d$ differs from $b_{2}%
(c_{1},\ldots,c_{k},d,\ldots,d)=g_{2}\circ d$. It remains to settle the case
$Q_{1}\neq Q_{2}$. Again, we can choose elements $c_{k+1},\ldots,c_{n}%
\in\mathbf{3}$ with $h_{1}=q_{1}(c_{k+1},\ldots,c_{n})\neq q_{2}%
(c_{k+1},\ldots,c_{n})=h_{2}$.

\noindent\textit{Case }17\textit{:} 0 and 2 are left zero elements, whence
$c_{k+1}=1$, and we can assume $h_{1}=0,\,h_{2}=1$. Now $b_{1}(1,\ldots
,1,c_{k+1},\ldots,c_{n})=1\circ0=0\neq1=1\circ1=b_{2}(1,\ldots,1,c_{k+1}%
,\ldots,c_{n})$.

\noindent\textit{Cases }16\textit{\ and }178\textit{:} for distinct elements
$h_{1},h_{2}\in\mathbf{3}$ there exists $e\in\mathbf{3}$ with $e\circ
h_{1}\neq e\circ h_{2}$. Hence it follows $b_{1}(e,\ldots,e,c_{k+1}%
,\ldots,c_{n})\neq b_{2}(e,\ldots,e,c_{k+1},\ldots,c_{n})$, concluding the proof.

\bigskip

\noindent5.2.\textit{\ The operation }33\textit{\ is Catalan. For
}35\textit{\ and }68\textit{, }$s(n)=2^{n-2}$\textit{.}

\smallskip

Consider a groupoid $(G,\circ)$ with idempotent elements $d,e(\neq d),f$ such that

\smallskip

\noindent$\left(  \alpha\right)  $ in the Cayley table of $\circ$, $d$ occurs
only in its own row;

\noindent$\left(  \beta\right)  $ in the row of $e$, $e\circ d$ occurs only once;

\noindent$\left(  \gamma\right)  $ $f$ is a right unit element.

\smallskip

First check that $\mathbf{3}$ with 33, 35, or 68 satisfies these conditions.
Now let $B_{1}=(P_{1}Q_{1})$ and $B_{2}=(P_{2}Q_{2})$ be bracketings of size
$n$ such that their induced operations over $\circ$ coincide. We prove
$p_{1}=p_{2}$. Suppose $k=|P_{1}|<|P_{2}|=l$. Then
\begin{align*}
b_{2}(e,\ldots,e,d,\ldots,d)  &  =p_{2}(e,\ldots,e)\circ q_{2}(d,\ldots
,d)=e\circ d,\\
b_{1}(e,\ldots,e,d,\ldots,d)  &  =p_{1}(e,\ldots,e)\circ q_{1}(e,\ldots
,e,d,\ldots,d).
\end{align*}
As by $\left(  \alpha\right)  $ we have $q_{1}(e,\ldots,e,d,\ldots,d)\neq d$,
from $\left(  \beta\right)  $ it follows that $b_{1}(e,\ldots,e,\allowbreak
d,\ldots,d)\neq b_{2}(e,\ldots,e,d,\ldots,d)$. Thus $|P_{1}|=|P_{2}|$, and ,
by $\left(  \gamma\right)  $, for arbitrary $c_{1},\ldots,c_{k}\in G$ it holds
$p_{1}(c_{1},\ldots,c_{k})=b_{1}(c_{1},\ldots,c_{k},f,\ldots,f)=b_{2}%
(c_{1},\ldots,c_{k},\allowbreak f,\ldots,f)=p_{2}(c_{1},\ldots,c_{k})$, which
was needed.

Take into account that 33 is surjective, and its Cayley table has no two
identical columns. We show that in the case of 33 if $b_{1}=b_{2}$, then
$q_{1}=q_{2}$, which together with $p_{1}=p_{2}$ implies \textit{via}
induction that 33 is Catalan. Indeed, suppose that, although $b_{1}=b_{2}$,
there exist $c_{k+1},\ldots,c_{n}\in\mathbf{3}$ such that $q_{1}%
(c_{k+1},\ldots,c_{n})\neq q_{2}(c_{k+1},\ldots,c_{n})$. Then the columns of
these two elements are also distinct, \textit{i.e.} $c\circ q_{1}%
(c_{k+1},\ldots,c_{n})\neq c\circ q_{2}(c_{k+1},\ldots,c_{n})$ for some
$c\in\mathbf{3}$. In virtue of 4.2.1 we can choose $c_{1},\ldots,c_{k}%
\in\mathbf{3}$ so that $p_{1}(c_{1},\ldots,c_{k})=c$. Now $b_{1}(c_{1}%
,\ldots,c_{n})=p_{1}(c_{1},\ldots,c_{k})\circ q_{1}(c_{k+1},\ldots,c_{n})\neq
p_{1}(c_{1},\ldots,c_{k})\circ q_{2}(c_{k+1},\ldots,c_{n})=b_{2}(c_{1}%
,\ldots,c_{n})$, a contradiction.

Concerning 35 and 68, observe that in these cases if $u\circ v\neq u\circ w$
then at least one of $v$ and $w$ is a left zero which satisfies $\left(
\alpha\right)  $. We have seen that $b_{1}=b_{2}$ implies $p_{1}=p_{2}$; now
we prove that the converse implication also holds. Suppose not, \textit{i.e.},
there are $c_{1},\ldots,c_{n}\in\mathbf{3}$ such that $b_{1}(c_{1}%
,\ldots,c_{n})=p_{1}(c_{1},\ldots,c_{k})\circ q_{1}(c_{k+1},\ldots,c_{n})\neq
p_{1}(c_{1},\ldots,c_{k})\circ q_{2}(c_{k+1},\ldots,c_{n})=b_{2}(c_{1}%
,\ldots,c_{n})$. Hence, without loss of generality, the element $d=q_{1}%
(c_{k+1},\ldots,c_{n})$ is a left zero, and $d$ does not occur in other rows.
We infer that $c_{k+1}=d$, and, as a consequence, $q_{2}(c_{k+1},\ldots
,c_{n})=d=q_{1}(c_{k+1},\ldots,c_{n})$, whence $b_{1}(c_{1},\ldots
,c_{n})=b_{2}(c_{1},\ldots,c_{n})$, a contradiction.

This shows that, for 35 and 68, $s(n)=s(n-1)+\cdots+s(2)+s(1)$, and this means
$s(n)=2^{n-2}$, as stated.

\bigskip

\noindent5.3.\textit{\ For the operation }624\textit{, }$s(n)=\lfloor
2^{n}/3\rfloor$\textit{.}

\smallskip

624 is actually $2x+2y\,\operatorname{mod}3$ on $\mathbf{3}$. We shall write
it in form $-x-y$; our considerations are valid for this operation on numbers,
too. An $n$-ary regular operation (over $-x-y$) is always of form
$t(x_{1},\ldots,x_{n})=\pm x_{1}\pm\cdots\pm x_{n}$. We call such operations
\textit{complete linear}. As $x_{1}-x_{2}+x_{3}$ shows, not every complete
linear operation is regular. Denote by $\pi(t)$ the number of $+$ signs in a
complete linear operation $t=t(x_{1},\ldots,x_{n})$, and call a complete
linear $t$ \textit{subregular}, if $\pi(t)\equiv2n-1\,(\operatorname{mod}3)$.
The following assertion can be checked immediately:

\smallskip

\noindent5.3.1. If $t,t_{1},t_{2}$ are complete linear operations such that
the equality $t(x_{1},\ldots,x_{n})=-t_{1}(x_{1},\ldots,x_{k})-t_{2}%
(x_{k+1},\ldots,x_{n})$ holds, then every one of $t,t_{1},t_{2} $ is
subregular provided the other two of them are subregular.

\smallskip

Next we characterize the regular operations over $-x-y$.

\smallskip

\noindent5.3.2. A complete linear operation $t(x_{1},\ldots,x_{n})$ is regular
over $-x-y$ if and only if it is subregular but not of form
\begin{equation}
x_{1}-x_{2}+x_{3}-\cdots+x_{n} \tag{5}%
\end{equation}
(\textit{i.e.}, not of odd arity with alternating signs and beginning with a
$+$ sign).

\smallskip

Clearly, this is true for $n\leq3$. Suppose that $t$ is regular. Then
$t(x_{1},\ldots,x_{n})=-t_{1}(x_{1},\ldots,x_{k})-t_{2}(x_{k+1},\ldots,x_{n})$
with $t_{1}$ and $t_{2}$ regular. By induction, $t_{1}$ and $t_{2}$ are
subregular, and 5.3.1 implies that $t$ is subregular. If $t$ is regular and it
is of form (5), then one of $t_{1}$ and $t_{2}$ --- say, $t_{1}$ --- must be
of even arity with alternating signs. However, a complete linear operation $t$
of arity $2m$ with alternating signs cannot be subregular, as $\pi
(t)=m\not \equiv 2\cdot2m-1\,(\operatorname{mod}3)$. Hence $t_{1}$ is not
subregular, a contradiction.

Conversely, assume that $t$ is subregular but not regular. We have to prove
that $t$ is of form (5). We show that the first sign in $t$ is $+$ . If not,
then $t(x_{1},\ldots,x_{n})=-x_{1}\pm x_{2}\pm\cdots\pm x_{n}=-x_{1}-(\mp
x_{2}\mp\cdots\mp x_{n})=-x_{1}-t_{2}(x_{2},\ldots,x_{n})$, and from 5.3.1 it
follows that $t_{2}$ is subregular. If, in addition, $t_{2}$ is not of form
(5), then by induction $t_{2}$ is regular, hence $t$ is regular, in contrary
to the assumption. However, if $t_{2}$ is of form (5), then%
\begin{align*}
t(x_{1},\ldots,x_{n})  &  =-x_{1}-x_{2}+x_{3}-\cdots+x_{n-1}-x_{n}=\\
&  =-(x_{1}+x_{2}-x_{3}+\cdots-x_{n-1})-x_{n}=\\
&  =-t_{1}(x_{1},\ldots,x_{n-1})-x_{n},
\end{align*}
and here $t_{1}$ is regular, implying again the regularity of $t$.

Thus, $t$ starts with a $+$ sign, and it is enough to prove that the signs
alternate in $t$. If not, consider the first two consecutive identical signs
in $t$. Suppose they are $+$ ; the other case can be treated analogously. Then%
\begin{align*}
t(x_{1},\ldots,x_{n})=  &  \ x_{1}-x_{2}+\cdots-x_{2k-2}+x_{2k-1}+x_{2k}\pm\\
&  \pm x_{2k+1}\pm\cdots\pm x_{n}=\\
=  &  -(-x_{1}+x_{2}-\cdots+x_{2k-2}-x_{2k-1}-x_{2k})-\\
&  -(\mp x_{2k+1}\mp\cdots\mp x_{n})=\\
=  &  -t_{1}(x_{1},\ldots,x_{2k})-t_{2}(x_{2k+1},\ldots,x_{n}).
\end{align*}
We can check that $t_{1}$ is subregular and not of form (5), hence regular;
further, $t_{2}$ is subregular by 5.3.1. As above, supposing that $t_{2}$ is
not of form (5) leads to a contradiction. Hence $t_{2}(x_{2k+1},\ldots
,x_{n})=x_{2k+1}-x_{2k+2}+\cdots-x_{n-1}+x_{n}$, and
\begin{align*}
t(x_{1},\ldots,x_{n})=  &  \ x_{1}-x_{2}+x_{3}-\cdots+x_{2k-1}+x_{2k}%
-x_{2k+1}+\\
&  +x_{2k+2}-\cdots+x_{n-1}-x_{n}=\\
=  &  -(-x_{1}+x_{2}-x_{3}+\cdots-x_{2k-1}-x_{2k}+x_{2k+1}-\\
&  -x_{2k+2}+\cdots-x_{n-1})-x_{n}=\\
=  &  -{t_{1}}^{\prime}(x_{1},\ldots,x_{n-1})-x_{n}.
\end{align*}
Here ${t_{1}}^{\prime}$ is subregular and not of form (5), so it is regular by
induction, whence we obtain that $t$ is regular, and this final contradiction
proves that a subregular but not regular complete linear operation is of form (5).

From 5.3.2 it follows that the number $s(n)$ of the $n$-ary regular operations
over $-x-y$ equals $\sum_{k}{\binom{n}{{3k+i}}}-(n\,\operatorname{mod}2)$, if
$n\equiv2-i\,(\operatorname{mod}3)~(i=0,1,2)$. It is known that each of these
numbers is equal to $\lfloor{2}^{n}/3\rfloor$ (see [6], Ch. 5, Exercise 75).
This completes the description of spectra of three-element groupoids with
minimal clones.

\smallskip

The next seven operations are of some interest from various reasons. The first
two pairs have the same spectra but with different coincidences of induced
regular operations. Fibonacci numbers appear at the fifth one. Nest structure
is exploited in the next example, and the last one is related to the Sheffer
operation on $\mathbf{2}$. These operations are numbered by their codes in the
Siena Catalog [2]:%

\begin{gather*}
\renewcommand{\arraystretch}{1.4}\setlength{\tabcolsep}{0.2cm}%
\begin{tabular}
[c]{cccc}%
$%
\begin{tabular}
[c]{|ccc}\hline
$0$ & $0$ & $2$\\
$0$ & $0$ & $2$\\
$2$ & $2$ & $1$%
\end{tabular}
$ & $%
\begin{tabular}
[c]{|ccc}\hline
$0$ & $0$ & $0$\\
$0$ & $0$ & $0$\\
$1$ & $0$ & $0$%
\end{tabular}
$ & $%
\begin{tabular}
[c]{|ccc}\hline
$0$ & $0$ & $1$\\
$0$ & $0$ & $1$\\
$1$ & $1$ & $0$%
\end{tabular}
$ & $%
\begin{tabular}
[c]{|ccc}\hline
$1$ & $1$ & $1$\\
$2$ & $2$ & $2$\\
$0$ & $0$ & $0$%
\end{tabular}
$\\
$1066$ & $10$ & $405$ & $3242$%
\end{tabular}
\\
\renewcommand{\arraystretch}{1.4}\setlength{\tabcolsep}{0.2cm}%
\begin{tabular}
[c]{ccc}%
$%
\begin{tabular}
[c]{|ccc}\hline
$0$ & $0$ & $0$\\
$0$ & $1$ & $0$\\
$0$ & $0$ & $1$%
\end{tabular}
$ & $%
\begin{tabular}
[c]{|ccc}\hline
$0$ & $0$ & $0$\\
$0$ & $1$ & $0$\\
$0$ & $1$ & $2$%
\end{tabular}
$ & $%
\begin{tabular}
[c]{|ccc}\hline
$1$ & $0$ & $0$\\
$0$ & $2$ & $0$\\
$0$ & $0$ & $0$%
\end{tabular}
$\\
$79$ & $82$ & $2407$%
\end{tabular}
\end{gather*}

\bigskip

\noindent5.4.\textit{\ For the operations }1066\textit{\ and }10\textit{,
}$s(n)=n-1$\textit{.}

\smallskip

Denote by $t(c_{1},\ldots,c_{n})$ the number of occurrences of $2$ among
$c_{1},\ldots,c_{n}$. Concerning 1066, induction shows that, for arbitrary
bracketing $B=(PQ)$ with $|B|=n,\,|P|=k$, and $c_{1},\ldots,c_{n}\in
\mathbf{3}$,
\[
b(c_{1},\ldots,c_{n})=2\text{ if and only if }t(c_{1},\ldots,c_{n})\text{ is
odd,}%
\]
and
\[
b(c_{1},\ldots,c_{n})=1\text{ if and only if both }t(c_{1},\ldots,c_{k})\text{
and }t(c_{k+1},\ldots,c_{n})\text{ are odd.}%
\]
As a consequence, $b(c_{1},\ldots,c_{n})=0$ if and only if both $t(c_{1}%
,\ldots,c_{k})$ and $\linebreak t(c_{k+1},\ldots,c_{n})$ are even. Hence it
follows that two bracketings of equal size induce the same operation if and
only if the sizes of their left factors are equal.

\smallskip

In order to manage 10 (which, for this once, will be written as
multiplication), we introduce the \textit{priority of a bracketing} $B$
($\mathrm{pr}(B)$ in sign) for $|B|>2$ as follows: If $B=(PQ)$ and $|P|>1$,
then $\mathrm{pr}(B)=0$; if $B=(x_{1}(x_{2}(\ldots(x_{k}(R))\ldots)))$, and
$\mathrm{pr}(R)=0$ or $|R|=2$, then $\mathrm{pr}(B)=k$. We call the bracketing
$R$ \textit{the core of} $B$. Clearly, if $n>2$, for every $k=0,1,\ldots,n-2$
there exist bracketings of size $n$ with priority $k$. Hence it is sufficient
to prove that two bracketings of size $n$ induce the same regular operation
over 10 if and only if they are of the same priority.

\textquotedblleft If\textquotedblright: $\mathrm{pr}(B)=0$ implies that $b$ is
the constant $0$ operation. If $k=n-2$ or $k=n-3$, then there is only one
bracketing $B$ with $\mathrm{pr}(B)=k$. Suppose $B_{1}$ and $B_{2}$ are of
size $n$ with cores $R_{1}$, resp. $R_{2}$, and $\mathrm{pr}(B_{1}%
)=\mathrm{pr}(B_{2})=k<n-3$. Then%
\begin{align*}
b_{1}(c_{1},\ldots,c_{n})  &  =(c_{1}(\ldots(c_{k}\cdot r_{1}(c_{k+1}%
,\ldots,c_{n}))\ldots))=(c_{1}(\ldots(c_{k}\cdot0)\ldots))=\\
&  =(c_{1}(\ldots(c_{k}\cdot r_{2}(c_{k+1},\ldots,c_{n}))\ldots))=b_{2}%
(c_{1},\ldots,c_{n})
\end{align*}
for arbitrary $c_{1},\ldots,c_{n}\in\mathbf{3}$.

\textquotedblleft Only if\textquotedblright: Let again $B_{1}$ and $B_{2}$ be
bracketings with cores as above, and let $\mathrm{pr}(B_{1})=k<l=\mathrm{pr}%
(B_{2})$. Induction on priority shows that bracketings with positive priority
induce nonconstant operations over 10. Hence there are $c_{k+1},\ldots
,c_{l},\allowbreak c_{l+1},\ldots,c_{n}\in\mathbf{3}$ such that $(c_{k+1}%
(\ldots(c_{l}\cdot r_{2}(c_{l+1},\ldots,c_{n}))\ldots))=1$. For $i=0,1$, check
the equality $(2(2(\ldots(2\cdot i)\ldots)))=(k-i)\,\operatorname{mod}2$,
where $k$ is the number of occurrences of 2 in the left side, and choose
$c_{1}=\cdots=c_{k}=2$. It follows%
\begin{align*}
b_{1}(c_{1},\ldots,c_{n})  &  =(c_{1}(\ldots(c_{k}\cdot r_{1}(c_{k+1}%
,\ldots,c_{n}))\ldots))=(c_{1}(\ldots(c_{k}\cdot0)\ldots))=\\
&  =k\,\operatorname{mod}2\neq(k-1)\,\operatorname{mod}2=(c_{1}(\ldots
(c_{k}\cdot1)\ldots))=\\
&  =(c_{1}(\ldots(c_{k}(c_{k+1}(\ldots(c_{l}\cdot r_{2}(c_{l+1},\ldots
,c_{n}))\ldots)))\ldots))=\\
&  =b_{2}(c_{1},\ldots,c_{n}).
\end{align*}

\bigskip

\noindent5.5.\textit{\ For the operations }405\textit{\ and }3242\textit{,
}$s(n)=3$\textit{\ if }$n>3$\textit{.}

\smallskip

Let $B_{1},B_{2}$ be bracketings of size $n$, $B_{i}=(P_{i}Q_{i})$. We show
that the induced regular operations $b_{1},b_{2}$ over 405 coincide if and
only if one of the following conditions is satisfied:

\smallskip

\noindent$\left(  1\right)  $ $|P_{1}|=|P_{2}|=1$;

\noindent$\left(  2\right)  $ $1<|P_{1}|,|P_{2}|<n-1$;

\noindent$\left(  3\right)  $ $|P_{1}|=|P_{2}|=n-1$.

\smallskip

Indeed, in the case $\left(  1\right)  $ the first variable, and in the case
$\left(  3\right)  $ the last variable determines the value of $b_{i}$. In the
case $\left(  2\right)  $ $b_{i}$ is the constant zero operation. Finally, if
$B_{1}=(x_{1}Q_{1}),\,B_{2}=(P_{2}x_{n})$, then $b_{1}(0,\ldots,2)=0\neq
1=b_{2}(0,\ldots,2)$.

\smallskip

3242 is $x+1\,\operatorname{mod}3$. Similarly to 4.1, for any bracketing $B$
and its induced operation $b$ over 3242 we have $b_{1}(c_{1},\ldots
,c_{n})=c_{1}+d\,\operatorname{mod}3$, where $d$ is the left depth of $x_{1}$
in $B$.

\bigskip

\noindent5.6.\textit{\ For the operation }79\textit{, }$s(n)=F_{n+1}%
-1$\textit{, where }$F_{k}$\textit{\ is the }$k\hspace{1pt}$\textit{th
Fibonacci number.}

\smallskip

First we show that, for bracketings $B_{1},B_{2}$ of equal size, $b_{1}$
coincides with $b_{2}$ if and only if the eggs of nests of $B_{1}$ are the
same as the eggs of nests of $B_{2}$. Suppose that $x_{i},x_{i+1}$ are the
eggs of a nest of $B_{1}$ but of no nest of $B_{2}$. Put $c_{j}=2$, if $j=i$
or $j=i+1$, and $c_{j}=1$ otherwise. Then $b_{1}(c_{1},\ldots,c_{n}%
)=1\neq0=b_{2}(c_{1},\ldots,c_{n})$. On the other hand, if the eggs of nests
of $B_{1}$ and $B_{2}$ are the same, induction on the number of nests proves
$b_{1}=b_{2}$. Note that this number is 1 exactly when $B_{1}$ and $B_{2}$ are
nests, and for nests we can apply the usual induction on size.

Choose several non-overlapping pairs $(i,i+1)$ in the sequence $1,\ldots,n$.
The number of such choices (including the empty choice) is $F_{n+1}$.
Induction shows that for every such nonempty choice $C$ there exists a
bracketing $B$ such that $x_{i},x_{i+1}$ are the eggs of a nest of $B$ if and
only if $(i,i+1)$ occurs in the choice $C$. This proves our proposition.

\bigskip

\noindent5.7.\textit{\ The operation }82\textit{\ is Catalan.}

\smallskip

Induction shows that the first (\textit{i.e.}, leftmost) right parenthesis in
$B$ together with its left pair encloses just the eggs of the leftmost
nontrivial (maximal) nest of $B$. Let $|B_{1}|=|B_{2}|=n,~b_{1}=b_{2}$, and
let the eggs in question of $B_{1}$ and $B_{2}$ consist of $x_{k},x_{k+1}$ and
$x_{l},x_{l+1}~(k<l)$, respectively. For $c_{1}=\cdots=c_{k}=c_{k+2}%
=\cdots=c_{n}=1,~c_{k+1}=2$ we get $b_{1}(c_{1},\ldots,c_{n})=0\neq
b_{2}(c_{1},\ldots,c_{n})$. Thus, the first right parentheses in $B_{1}$ and
$B_{2}$ cannot be in different positions. Collapsing $x_{k}$ and $x_{k+1}$ we
obtain quotient bracketings $B_{1}^{\prime}$ and $B_{2}^{\prime}$ of size
$n-1$. Remark that, for arbitrary $c_{1},\ldots,c_{k-1},c_{k+1},\ldots
,c_{n}\in\mathbf{3}$,$~b_{i}^{\prime}(c_{1},\ldots,c_{k-1},c_{k+1}%
,\ldots,c_{n})=b_{i}(c_{1},\ldots,c_{k-1},2,c_{k+1},\ldots,c_{n})$ holds, as
$2$ is a left unit for 82. In such a way, $b_{i}$ determines $b_{i}^{\prime}$,
and the latter determines the place of the first right parenthesis in
$B_{i}^{\prime}$, which is the second right parenthesis in $B_{i}$;
\textit{etc}. We see that the induced operation determines the positions of
all right parentheses in its parent bracketing. Now 5.7 follows from 2.6.

\bigskip

\noindent5.8.\textit{\ The operation }2407\textit{\ is Catalan.}

\smallskip

The proof consists of a suitable adaptation of 4.2. The observations 4.2.1,
4.2.2, and 4.2.3 apply to 2407. Now, from $B_{1}=(P_{1}Q_{1}),~B_{2}%
=(P_{2}Q_{2})$, and $b_{1}=b_{2}$ we can deduce not only the equivalence of
$p_{1}(c_{1},\ldots,c_{k})=0$ and $p_{2}(c_{1},\ldots,c_{l})=0$ but also that
of $p_{1}(c_{1},\ldots,c_{k})=1$ and $p_{2}(c_{1},\ldots,c_{l})=1$. Thus,
again we have $p_{1}=p_{2}$, and, by induction, $P_{1}=P_{2}$. In order to
refute $Q_{1}\neq Q_{2}$, assume that there exist $c_{k+1},\ldots,c_{n}%
\in\mathbf{3}$ with $q_{1}(c_{k+1},\ldots,c_{n})=i\neq j=q_{2}(c_{k+1}%
,\ldots,c_{n})$; here we can suppose $i\neq2$. There are $c_{1},\ldots
,c_{k}\in\mathbf{3}$ with $p_{1}(c_{1},\ldots,c_{k})=i$. Then $b_{1}%
(c_{1},\ldots,c_{n})=i\circ i=i+1\,\operatorname{mod}3\neq i\circ
j=b_{2}(c_{1},\ldots,c_{n})$.

\smallskip

The Sheffer function on $\mathbf{2}$ and 2407 on $\mathbf{3}$ are the smallest
instances of groupoids $(\mathbf{n},\circ)$ with operations%
\begin{equation}
i\circ j=%
\begin{cases}
i+1, & \text{if }i=j\\
0, & \text{otherwise.}%
\end{cases}
\tag{6}%
\end{equation}
All these groupoids are \textit{primal}\thinspace; \textit{i.e.}, all possible
operations on $\mathbf{n}$ are term operations of such a groupoid. The proof
of 5.8 can be generalized for them without trouble. Hence we could (in fact,
we did) conjecture for a minute that primality implies a Catalan spectrum;
however, operation 3233 testifies that this is not the case. Its Cayley table
comes from that of 3242 by writing $1\circ2=0$ instead of $1\circ2=2$. For
3233 we have $s_{6}=41<C_{5}(=42)$. Actually,%
\[
x_{1}\circ\left(  \left(  x_{2}\circ\left(  x_{3}\circ\left(  x_{4}\circ
x_{5}\right)  \right)  \right)  \circ x_{6}\right)  =x_{1}\circ\left(  \left(
x_{2}\circ\left(  \left(  x_{3}\circ x_{4}\right)  \circ x_{5}\right)
\right)  \circ x_{6}\right)
\]
identically holds for 3233 on $\mathbf{3}$ (but no other regular operations
over 3233 induced by distinct bracketings of size $\leq6$ are equal). On the
other hand, the primality of $\mathbf{3}$ with 3233 as well as of $\mathbf{n}$
with operation (6) follows, \textit{e.g.}, from Rousseau's criterion: a finite
algebra with a single operation is primal if and only if it has neither proper
subalgebras, nor congruences, nor automorphisms [12].

We have checked all the 3330 entries of the Siena Catalog by computer for the
five initial elements of their spectra, \textit{i.e.}
$(s(3),s(4),s(5),s(6),s(7))$. It is known that there are 24 nonisomorphic
three-element semigroups. The table below shows the number of essentially
distinct three-element nonassociative groupoids with a given initial segment
of spectrum:%

\[
\renewcommand{\arraystretch}{1.2} \setlength{\tabcolsep}{0.2cm}%
\begin{tabular}
[c]{rrrrrrr}%
$2$ & $2$ & $2$ & $2$ & $2$ &  & $16$\\
$2$ & $3$ & $3$ & $3$ & $3$ &  & $4$\\
$2$ & $3$ & $4$ & $5$ & $6$ &  & $15$\\
$2$ & $4$ & $4$ & $4$ & $4$ &  & $2$\\
$2$ & $4$ & $5$ & $6$ & $7$ &  & $6$\\
$2$ & $4$ & $6$ & $8$ & $10$ &  & $4$\\
$2$ & $4$ & $7$ & $12$ & $20$ &  & $4$\\
$2$ & $4$ & $7$ & $12$ & $21$ &  & $12$\\
$2$ & $4$ & $8$ & $15$ & $27$ &  & $12$\\
$2$ & $4$ & $8$ & $16$ & $32$ &  & $62$\\
$2$ & $5$ & $8$ & $12$ & $16$ &  & $2$\\
$2$ & $5$ & $10$ & $18$ & $31$ &  & $4$\\
$2$ & $5$ & $10$ & $20$ & $40$ &  & $4$%
\end{tabular}
\qquad\qquad\qquad%
\begin{tabular}
[c]{rrrrrrr}%
$2$ & $5$ & $10$ & $21$ & $42$ &  & $5$\\
$2$ & $5$ & $11$ & $23$ & $47$ &  & $2$\\
$2$ & $5$ & $11$ & $24$ & $53$ &  & $4$\\
$2$ & $5$ & $12$ & $28$ & $65$ &  & $12$\\
$2$ & $5$ & $13$ & $34$ & $87$ &  & $12$\\
$2$ & $5$ & $13$ & $34$ & $89$ &  & $2$\\
$2$ & $5$ & $13$ & $34$ & $90$ &  & $4$\\
$2$ & $5$ & $13$ & $34$ & $91$ &  & $24$\\
$2$ & $5$ & $13$ & $35$ & $96$ &  & $2$\\
$2$ & $5$ & $13$ & $35$ & $97$ &  & $32$\\
$2$ & $5$ & $14$ & $41$ & $123$ &  & $6$\\
$2$ & $5$ & $14$ & $41$ & $124$ &  & $16$\\
$2$ & $5$ & $14$ & $42$ & $132$ &  & $3038$%
\end{tabular}
\]

\noindent Several sequences beginning with some quintuples above,
\textit{e.g.} $(2,5,10,21,42)$ (\textit{cf.} 5.3) and $(2,5,14,41,123)$, are
recently missing in the Encyclopedia [13].

\section{General remarks and problems}

\smallskip

All the spectra considered up to now are monotonic. Groups with the commutator
operation provide examples of nonmonotonic spectra: if a group $G$ is
nilpotent then there exists an $n$ such that all $n$-ary regular term
operations over the commutator of $G$ are equal (to the constant unit
operation), hence $s(n)=1$, and if $G$ is not nilpotent of class 2 then the
commutator is not associative (see, \textit{e.g.} [10]
%, Vol. 2., p. 179, and Vol. 1., pp. 99---102, respectively
). The spectrum always stabilizes in these examples: $s(n)=1$ implies $s(m)=1
$ for every $m>n$. In fact, this is a common property of all spectra, which
generalizes the generalized associative law:

\bigskip

\noindent6.1.\textit{\ For an arbitrary spectrum }$s$\textit{,~}%
$s(n)=1$\textit{\ for some }$n\,(\geq3)$\textit{\ implies }$s(m)=1$%
\textit{\ for every }$m>n$\textit{.}

\smallskip

Call two bracketings of size $m$ \textit{adjacent} if there exists a $j$ such
that $x_{j},x_{j+1}$ are eggs of nests for each of these bracketings. It is
easy to see that the transitive closure of the adjacency relation is the
trivial equivalence if $m\geq5$.

Let $n\,(\geq3)$ be a number such that $s(n)=1$ for an operation $\circ$ on a
set $M$. Consider bracketings $B,B^{\ast}$ of size $n+1$. We have to prove
$b=b^{\ast}$. For $n=3$ this is the generalized associative law. Assume $n>3$.
Then $n+1\geq5$, hence there exist bracketings $B_{0}=B,B_{1},\ldots
,B_{k}=B^{\ast}$ such that, for $i=0,1,\ldots,k-1$,$~B_{i}$ is adjacent to
$B_{i+1}$. Let $x_{j},x_{j+1}$ be common eggs of a nest of $B_{i}$ and a nest
of $B_{i+1}$. Replacing $(x_{j}x_{j+1})$ by $x_{j} $ in both of them, we
obtain quotient bracketings $B_{i}^{\prime},B_{i+1}^{\prime}$ of size $n$. As
$s(n)=1$, we have $b_{i}^{\prime}=b_{i+1}^{\prime}$, and thus%
\begin{align*}
b_{i}(c_{1},\ldots,c_{n+1})  &  =b_{i}^{\prime}(c_{1},\ldots,c_{j-1}%
,c_{j}\circ c_{j+1},c_{j+2},\ldots,c_{n+1})=\\
&  =b_{i+1}^{\prime}(c_{1},\ldots,c_{j-1},c_{j}\circ c_{j+1},c_{j+2}%
,\ldots,c_{n+1})=\\
&  =b_{i+1}(c_{1},\ldots,c_{n+1})
\end{align*}
for arbitrary $c_{1},\ldots,c_{n+1}\in M$.

\smallskip

Groups provide also examples showing that the difference $s(n)-s(n-1)$ of
consecutive entries of a spectrum can be arbitrarily large:

\bigskip

\noindent6.2\textit{. The spectrum of the commutator operation on the dihedral
group of degree }$2^{t}~(t\geq3)$\textit{\ is}%
\[
s(n)=%
\begin{cases}
2, & \text{if }n=3\\
n, & \text{if }3<n\leq t\\
1, & \text{if }n>t.
\end{cases}
\]

$D_{m}$, the dihedral group of degree $m$ is generated by a rotation $\alpha$
of order $m$ and a reflection $\rho$. We write $i$ for $\alpha^{i}$ and
$j^{\prime}$ for $\alpha^{j}\rho$. Here is the concise Cayley table of the
commutator on $D_{m}$:

\smallskip%

\[
\renewcommand{\arraystretch}{1.4}\setlength{\tabcolsep}{0.2cm}%
\begin{tabular}
[t]{c|c|c}
& $j$ & $j^{\prime}$\\\hline
$i$ & $0$ & $-2i\,\operatorname{mod}m$\\\hline
$i^{\prime}$ & $2j\,\operatorname{mod}m$ & $2\left(  i-j\right)
\,\operatorname{mod}m$%
\end{tabular}
\]

\bigskip

The following observations are immediate: If a bracketing $B$ over the
commutator on $D_{n}$ has at least two nests, then it induces the constant
zero operation. Further, if $B$ is a nest with eggs $x_{k},\,x_{k+1}$, then
$b(c_{1},\ldots,c_{n})\neq0$ only if all $c_{i}\,(\in D_{m})$ but at most one
of $c_{k},c_{k+1}$ are of form $i^{\prime}$ (\textit{i.e.,} $\alpha^{i}\rho$).
From the Cayley table we learn that for such a nest $B$ and such elements
$c_{1},\ldots,c_{n}$%
\begin{equation}
b(c_{1},\ldots,c_{n})=[c_{k},c_{k+1}]\,{2^{k-1}}(-2)^{n-k-1}%
\,\operatorname{mod}m \tag{7}%
\end{equation}
holds. From (7) we infer that the position of eggs of $B$ determines the
induced operation $b$. As all commutators are of form $2u\,\operatorname{mod}%
m$, (7) shows also that always $b(c_{1},\ldots,c_{n})={2^{n-1}}\cdot
v\,\operatorname{mod}m$ with suitable integers $v$. This means that $b$ is the
zero operation if $m=2^{t}$ and $n>t$.

It remains to show that nests of equal size $n\,(\leq t)$ but with distinct
eggs induce distinct operations. In fact, besides $B$ consider another nest
$B^{\prime}$ with eggs $x_{l},x_{l+1}~(l>k)$. Let $c_{k}=1,\,c_{k+1}%
=2^{\prime}$, and choose elements $c_{i}~(i\neq k,k+1)$ of form $i^{\prime}$
arbitrarily. Then $[1,2^{\prime}]=-2\,\operatorname{mod}2^{t}$, and, by (7),
$b(c_{1},\ldots,c_{n})=(-1)^{n-k}2^{n-1}\,\operatorname{mod}2^{t}\neq0$. On
the other hand, $l>k$ implies $b^{\prime}(c_{1},\ldots,c_{n})=0$ because
$c_{k}=1$, and $x_{k}$ is out of the egg of $B^{\prime}$.

The same reasoning shows that the commutator on $D_{1},\,D_{2}$ and $D_{4}$ is
associative, and if $m$ is not a power of $2$ (\textit{e.g.}, in the case of
$D_{3}=S_{3}$) the spectrum of the commutator on $D_{m}$ is $s(n)=n$ for $n>3$.

\medskip

The next example leads to groupoids whose spectra begin with arbitrarily many
Catalan numbers and still reach 1.

\bigskip

\noindent6.3.\textit{\ The following operation on the nonnegative integers is
Catalan:}%
\[
a\circ b=%
\begin{cases}
\min(a,b)-1, & \text{if }a,b>0\\
0, & \text{otherwise.}%
\end{cases}
\]

For the proof, denote by $d_{{}_{B}}(x_{i})$ the depth of $x_{i}$ in the
bracketing $B$. Consider an arbitrary bracketing $B=(PQ)$ with $|B|=n,~|P|=k$.
First we show that%
\[
b(d_{{}_{B}}(x_{1})+1,\ldots,d_{{}_{B}}(x_{n})+1)=1.
\]
Note that, for any $B$, $b(c_{1},\ldots,c_{n})>0$ implies $b(c_{1}%
+1,\ldots,c_{n}+1)=b(c_{1},\ldots,c_{n})+1$. By induction we have $p(d_{{}%
_{B}}(x_{1}),\ldots,d_{{}_{B}}(x_{k}))=p(d_{{}_{P}}(x_{1})+1,\ldots,d_{{}_{P}%
}(x_{k})+1)=1$, and similarly $q(d_{{}_{B}}(x_{k+1}),\ldots,d_{{}_{B}}%
(x_{n}))=1$, whence it follows%
\begin{align*}
b(d_{{}_{B}}(x_{1})+1,\ldots,d_{{}_{B}}(x_{n})+1)=  &  \ p(d_{{}_{B}}%
(x_{1})+1,\ldots,d_{{}_{B}}(x_{k})+1)\circ\\
&  \circ q(d_{{}_{B}}(x_{k+1})+1,\ldots,d_{{}_{B}}(x_{n})+1)=\\
=  &  \ (1+1)\circ(1+1)=1,
\end{align*}
as needed.

Next we show that for any other $B^{\prime}$ of size $n$, $b^{\prime}%
(d_{{}_{B}}(x_{1})+1,\ldots,d_{{}_{B}}(x_{n})+1)=0$. Again, induction shows
that for arbitrary $B$, nonnegative integers $c_{1},\ldots,c_{n}$, and
$i~(1\leq i\leq n)$%
\begin{equation}
b(c_{1},\ldots,c_{n})\leq\max(c_{i}-d_{{}_{B}}(x_{i}),0) \tag{8}%
\end{equation}
holds; we omit the details. As $B^{\prime}\neq B$, 2.7 implies that there
exists an $i$ such that $d_{{}_{B^{\prime}}}(x_{i})\neq d_{{}_{B}}(x_{i})$,
and in view of (1) we can suppose even $d_{{}_{B^{\prime}}}(x_{i})>d_{{}_{B}%
}(x_{i})$. Then applying (8) to $B^{\prime}$ we obtain%
\[
b^{\prime}(d_{{}_{B}}(x_{1})+1,\ldots,d_{{}_{B}}(x_{n})+1)\leq\max(d_{{}_{B}%
}(x_{i})+1-d_{{}_{B^{\prime}}}(x_{i}),0)=0,
\]
concluding the proof.

For any bracketing $B$ with $|B|=k<n$, and for every $i\,(=1,\ldots,k)$, we
have $d_{{}_{B}}(x_{i})<k$, hence $d_{{}_{B}}(x_{i})+1\in\mathbf{n}$.
Therefore the above reasoning shows that in $(\mathbf{n},\circ)$, which is a
subgroupoid of $(\mathcal{N}_{0},\circ)$, distinct bracketings of size
$k\,(<n)$ induce different regular operations. On the other hand, every
bracketing $B$ whose size exceeds $2^{n-2}$ has a symbol $x_{j}$ with
$d_{{}_{B}}(x_{j})\geq n-1$. Applying (8) to the regular operation $b$ of
$(\mathbf{n},\circ)$ we obtain%
\[
b(c_{1},\ldots,c_{n})\leq\max(c_{j}-d_{{}_{B}}(x_{j}),0)=0,
\]
as $c_{j}\leq n-1$. Hence any bracketing of size $2^{n-2}+1$ induces the
constant zero operation of $(\mathbf{n},\circ)$. Thus, for the spectrum of
$(\mathbf{n},\circ)$,~\thinspace$s(k)=C_{k-1}$ if $k<n$, and $s(k)=1$ if
$k>2^{n-2}$.

\smallskip

The study of spectra of linear operations $px+qy$ (and $px+qy+r$) on numbers
(or, more generally, on modules over rings) also offers remarkable facts. As a
specimen, we prove the following generalization of 3.2.

\bigskip

\noindent6.4.\textit{\ The linear operations }$px+py$\textit{\ and }%
$x+py$\textit{\ on the complex numbers are not Catalan if and only if }%
$p$\textit{\ is a root of unity.}

\smallskip

Concerning $px+py$, induction shows that for any bracketing $B$ of size $n$,
the induced operation over $px+py$ is%
\begin{equation}
b(x_{1},\ldots,x_{n})=\sum_{i=1}^{n}p^{d_{i}}x_{i}, \tag{9}%
\end{equation}
where $d_{i}$ is the depth of $x_{i}$ in $B$. From 2.7 it follows that if $p$
is not a root of unity then $px+py$ is Catalan. Suppose $p^{k}=1$. Define the
bracketings $B_{i}$ by $B_{1}=(xx)$, and $B_{n+1}=(B_{n}B_{n})$ for $n>0$. The
depth sequences of $B^{\prime}=(xB_{k})$ and $B^{\prime\prime}=(B_{k}x)$ are
$(1,k+1,\ldots,k+1)$ and $(k+1,\ldots,k+1,1)$, respectively. Now (9) implies
$b^{\prime}=b^{\prime\prime}$. Hence, for $m=2^{k}+1$, $s(m)<C_{m-1}$.

Analogous considerations apply to $x+py$: (9) remains valid for this case with
\textit{right} depths instead of depths. If $p$ is not a root of unity, 2.8
guarantees that $x+py$ is Catalan. Suppose again $p^{k}=1$, and redefine
$B_{i}$ by $B_{1}=(xx)$, and $B_{n+1}=(xB_{n})$ for $n>0$. The RD-sequences of
$B^{\prime}=(B_{k}x)$ and $B^{\prime\prime}=B_{k+1}$ are $(0,1,2,\ldots,k,1)$
and $(0,1,2,\ldots,k,k+1)$, respectively, implying $b^{\prime}=b^{\prime
\prime}$, and, for $m=k+2$, \thinspace$s(m)<C_{m-1}$.

\medskip

In conclusion, we formulate a few problems:

\smallskip

\noindent1. For every positive integer $n$ there exists a minimal $f(n)$ with
the property that, if for two spectra $s_{1},s_{2}$ of $n$-element groupoids
$s_{1}(i)=s_{2}(i)$ holds whenever $i\leq f(n)$, then these spectra coincide.
Propositions 4.1---3 imply $f(2)=4$, and the table at the end of Section 5
shows that $f(3)\geq7$. What is the actual value of $f(3)$ (and that of
$f(4)$, \textit{etc.})?

\smallskip

\noindent2. We gave a rough estimation for the subsequent entries of a
spectrum with a given initial segment in 2.3 which \textit{e.g.,} for $s(3)=2
$ and $s(4)=4$ provides $s(5)\leq12$. However, a case-by-case analysis shows
that $s(3)=2$ and $s(4)=4$ actually imply $s(5)\leq8$. Do they imply
$s(n)\leq2^{n-2}$ for all $n\,(>1)$\thinspace? If so, call $s(n)=2^{n-2}$ a
\textit{maximal extension} of the initial segment $(2,4)$. Prove or disprove
that the maximal extension of $(2,3)$ is $s(n)=n-1$, and that of $(2,2)$ is
$s(n)=2$.

\smallskip

\noindent3. All nonconstant spectra we exhibited above are ultimately constant
or monotonic. In the latter case their growth rates are either linear or
exponential. Is there any other possibility? More concretely: find,
\textit{e.g.}, a spectrum with quadratic growth rate.

\smallskip

\noindent4. The statistics of the three-element groupoids and the abundance of
appropriate examples leave such an impression that a huge majority of binary
operations is Catalan. Is it true that, in some sense, almost all operations
are Catalan (or almost Catalan)?


\begin{thebibliography}{99}                                                                                               %


\bibitem {b1}M. K. Bennett, G. Birkhoff, \textit{Two families of Newman
lattices}, Algebra Universalis, \textbf{32 }(1994), 115--144.

\bibitem {b2}J. Berman, S. Burris, \textit{A computer study of }%
$3$\textit{-element groupoids}, in: Logic and Algebra (Pontignano, 1994),
Lecture Notes in Pure and Appl. Math., 180, Dekker, 1996. (pp. 379-429)

\bibitem {b3}P. M. Cohn, \textit{Universal Algebra}, Harper \& Row, 1965.

\bibitem {b4}B. Cs\'{a}k\'{a}ny, \textit{Three-element groupoids with minimal
clones}, Acta Sci. Math. (Szeged), \textbf{45 }(1983), 111--117.

\bibitem {b5}F. G\"{o}bel, R. P. Nederpelt, \textit{The number of numerical
outcomes of iterated powers}, Amer. Math. Monthly, \textbf{78 }(1971), 1097--1103.

\bibitem {b6}R. L. Graham, D. E. Knuth, O. Patashnik, \textit{Concrete
Mathematics}, Addison-Wesley, 1994.

\bibitem {b7}R. K. Guy, J. L. Selfridge, \textit{The nesting and roosting
habits of laddered parentheses}, Amer. Math. Monthly, \textbf{80 }(1973), 868--876.

\bibitem {b8}N. Jacobson, \textit{Lectures in Abstract Algebra, Vol. I.}, D.
Van Nostrand, 1951.

\bibitem {b9}S. C. Kleene, \textit{Introduction to Metamathematics}, D. Van
Nostrand, 1952.

\bibitem {b10}A. G. Kurosh, \textit{The Theory of Groups, Vol. 1--2.},
Chelsea, 1955.

\bibitem {b11}L. Lov\'{a}sz, \textit{Combinatorial Problems and Exercises},
(2nd edition), North-Holland, 1993.

\bibitem {b12}G. Rousseau, \textit{Completeness in finite algebras with a
single operation}, Proc. Amer. Math. Soc., \textbf{18 }(1967), 1009--1013.

\bibitem {b13}N. J. A. Sloane, Simon Plouffe, \textit{The Encyclopedia of
Integer Sequences}, Academic Press, 1995.

\bibitem {b14}R. P. Stanley, \textit{Enumerative Combinatorics}, Vol. 2.,
Cambridge University Press, 1999.

\bibitem {b15}D. Tamari, \textit{The algebra of bracketings and their
enumeration}, Nieuw Arch. Wisk. (3), \textbf{10 }(1962), 131--146.
\end{thebibliography}
\end{document}